\documentclass[12pt]{article}
\usepackage{amsmath}
\usepackage{amssymb}
\usepackage[ctr1,eng]{myarticl}

\newcommand{\Aut }{\mathrm{Aut}}
\newcommand{\bigW }{\widetilde{W}}

\newcommand{\cB }{\mathcal{B}}

\newcommand{\End }{\mathrm{End}}

\newcommand{\lact }{.}

\newcommand{\lcoa }{\delta }

\newcommand{\Ndbasis }{\boldsymbol{\mathrm{e}}}

\newcommand{\ndN }{\mathbb{N}}
\newcommand{\ndQ }{\mathbb{Q}}
\newcommand{\ndR }{\mathbb{R}}
\newcommand{\ndZ }{\mathbb{Z}}

\newcommand{\ot }{\otimes }

\newcommand{\PBW }{Poincar\'e--Birkhoff--Witt }
\newcommand{\proots }{\boldsymbol{\Delta }^+}

\newcommand{\Rnvec }[1]{\boldsymbol{\mathrm{#1}}}

\newcommand{\roots }{\boldsymbol{\Delta }}
\newcommand{\scalp }[2]{\langle #1,#2\rangle _{\ndR ^n}}
\newcommand{\scp }[2]{p#1_{#2}}
\newcommand{\scq }[3]{{q#1}_{#2}^{(#3)}}
\newcommand{\SLZ }[1]{\mathrm{SL}(#1,\ndZ )}

\newcommand{\tr }{\mathrm{tr}}
\newcommand{\YD }{Yetter--Drinfel'd }

\title{Rank 2 Nichols algebras with finite arithmetic root system
\thanks{Supported by the European Community under a Marie Curie
Intra-European Fellowship}}
\author{I.~Heckenberger}

\begin{document}

\maketitle

\begin{abstract}
The concept of arithmetic root systems is introduced. It is shown
that there is a one-to-one correspondence between arithmetic root
systems and Nichols algebras of diagonal type having a finite set
of (restricted) \PBW generators. This has strong consequences for
both objects. As an application all rank 2 Nichols algebras of
diagonal type having a finite set of (restricted) \PBW generators
are determined.


\end{abstract}

\section{Introduction}

Nichols algebras play an essential role in the classification of
pointed Hopf algebras by the method of Andruskiewitsch and
Schneider \cite{a-AndrSchn98}. For this reason it is of particular
interest to determine all Nichols algebras with certain finiteness
conditions. Kharchenko \cite{a-Khar99} proved important results
about the general structure of such algebras. Further, a close
relation between Nichols algebras of diagonal type and the theory
of semi-simple Lie algebras was pointed out (see for example
\cite{a-AndrSchn00}). The latter are themselves characterized by
root systems. In this paper the concept of arithmetic root systems
is introduced, using mainly the Weyl--Brandt groupoid
\cite{a-Heck04c} associated to a bicharacter on $\ndZ ^n$. This
allows us to make the above relation precise.

Among Nichols algebras those of diagonal type are best understood.
There exist classification results for a subclass, namely Nichols
algebras of Cartan type \cite{a-Rosso98}, \cite{a-AndrSchn00},
\cite{a-Heck04c}. More general (published) examples of diagonal
type are known only in the rank 2 case, see \cite{inp-AndrSchn02}
and \cite{a-Heck04a}. One of the main aims of this paper is to
generalize the classification result for finite dimensional rank 2
Nichols algebras of diagonal type proved in \cite{a-Heck04b}, by
allowing a finite set of \PBW generators, but not requiring, that
all of them are restricted.

The paper is organized as follows. In Section \ref{sec-Ars} the
definition of arithmetic root systems is formulated. In Section
\ref{sec-WBg} the Weyl--Brandt groupoid of a Nichols algebra is
recalled. With Theorem \ref{t-corr} a one-to-one correspondence
between arithmetic root systems and Nichols algebras of diagonal
type with a finite set of (restricted) \PBW generators is given,
and two important consequences are noted. Finally, in Section
\ref{sec-Ars2} the second main assertion of the present paper is
formulated. Theorem \ref{t-class2} gives a list of all rank 2
Nichols algebras of diagonal type, which have a finite set of
(restricted) \PBW generators. The proof was possible mainly due to
the concept of Weyl equivalence of Nichols algebras introduced in
\cite{a-Heck04d}, and some technical results in
\cite[Sect.~4]{a-Heck04b}.

Throughout this paper $k$ denotes a field of characteristic zero
and tensor products $\ot $ are taken over this field. The set of
natural numbers including respectively not including zero is
denoted by $\ndN _0$ and $\ndN $, respectively.

I would like to thank the algebra and representation theory group
at the Weizmann Institute of Science in Rehovot, especially
A.~Joseph, A.~Melnikov, B.~Noyvert, and M.~Gorelik, for their
hospitality, their interest in the present research subject, and
numerous helpful discussions during my fellowship between February
2004 and January 2005.

\section{Arithmetic root systems}
\label{sec-Ars}

For a fixed $n\in \ndN $ let $\bigW $ denote the Brandt groupoid
\cite[Sect.~3.3]{b-ClifPres61}, \cite[Sect.~5]{a-Heck04c}
consisting of all pairs $(T,B)$ where $T\in \Aut (\ndZ ^n)$ and
$B$ is an ordered basis of $\ndZ ^n$, and the composition
$(T_1,B_1)\circ (T_2,B_2)$ is defined (and is then equal to
$(T_1T_2,B_2)$) if and only if $T_2(B_2)=B_1$. The Brandt groupoid
$\bigW $ is naturally acting on the set of all (ordered) bases of
$\ndZ ^n$ via the rule
\begin{align}\label{eq-WBaction}
(T,B)(B')=\begin{cases}
 T(B) & \text{if $B=B'$,}\\
 \text{not defined} & \text{otherwise.}
\end{cases}
\end{align}
Let $G$ be an abelian group and $\chi :\ndZ ^n\times \ndZ ^n\to G$
a bicharacter, i.\,e.~the map $\chi $ satisfies the properties
\begin{align*}
 \chi (0,e)=&\chi (e,0)=1,\\
 \chi (e'+e'',e)=&\chi (e',e)\chi (e'',e),\\
 \chi (e,e'+e'')=&\chi (e,e')\chi (e,e'')
\end{align*}
for all $e,e',e''\in \ndZ ^n$. Assume that $E=\{e_1,\ldots ,e_n\}$
is a basis of $\ndZ ^n$ and for an $i\in \{1,2,\ldots ,n\}$ and
all $j\in \{1,2,\ldots ,n\}$ with $j\not=i$ the numbers
\begin{align}\label{eq-mij}
\begin{aligned}
 m_{ij}:=\min \{m\in \ndN _0\,|\,&\text{either }\chi (e_i,e_i)^m\chi
 (e_i,e_j)\chi (e_j,e_i)=1\\
 &\text{or }\chi (e_i,e_i)^{m+1}=1,\chi (e_i,e_i)\not=1\}
\end{aligned}
\end{align}
exist. Then let $s_{i,E}\in \Aut (\ndZ ^n)$ denote the linear map
defined by
\begin{align}\label{eq-si}
s_{i,E}(e_j):=&
 \begin{cases}
 -e_i & \text{if $j=i$,}\\
 e_j+m_{ij}e_i & \text{if $j\not=i$.}
 \end{cases}
\end{align}
The map $s_{i,E}$ is a pseudo-reflection
\cite[Ch.~5,\S2]{b-BourLie4-6}, that is $\mathrm{rk}\,(s_{i,E}-\id
)=1$. Moreover it satisfies the equations
$s_{i,s_{i,E}(E)}=s_{i,E}$ and $s_{i,E}^2=\id $. Note that
$s_{i,E}$ doesn't depend on the antisymmetric part of $\chi $.

Let $E$ be an ordered basis of $\ndZ ^n$. Define $W_{\chi ,E}$ as
the smallest Brandt subgroupoid of $\bigW $ which contains $(\id
,E)$, and if $(\id, E')\in W_{\chi ,E}$ for an ordered basis $E'$
of $\ndZ ^n$ then $(s_{i,E'},E'),(\id ,s_{i,E'}(E'))\in W_{\chi
,E}$ whenever $s_{i,E'}$ is defined.
\begin{defin}\label{d-ars}
 Let $E$ be an ordered basis of $\ndZ ^n$ and $\chi $ a bicharacter on
 $\ndZ ^n$ such that $W_{\chi ,E}$ is finite. Assume that for all
 ordered bases $E'$ of $\ndZ ^n$ with $(\id ,E')\in W_{\chi ,E}$ the maps
 $s_{i,E'}\in \Aut (\ndZ ^n)$ are well-defined for all $i\in
 \{1,2,\ldots ,n\}$. Set $\roots :=\bigcup \{E'\,|\,(\id,E')\in
 W_{\chi ,E}\}\subset \ndZ ^n$ (union of sets, i.\,e.\ all
 elements appear with multiplicity one). The triple $(\roots ,\chi
 ,E)$ is called an \textit{arithmetic root system}.
\end{defin}
Note that $(\roots ,\chi ,E')=(\roots ,\chi ,E)$ for all bases
$E'$ of $\ndZ ^n$ such that $(\id ,E')\in W_{\chi ,E}$. Moreover,
any arithmetic root system has the following properties.
\begin{itemize}
 \item[(ARS1)] If $\alpha \in \roots $ then $\lambda \alpha \in
 \roots $ for some $\lambda \in \ndQ $ if and only if $\lambda ^2=1$.
\item[(ARS2)] If $G=k\setminus \{0\}$ for a field $k$ then one has
 $\roots =\proots _E\cup -\proots _E$, where $\proots _E=\{\alpha
 \in \roots \,|\,\alpha =\sum _{i=1}^nm_ie_i,\,m_i\in \ndN _0\text{
 for all }i\}$.
\end{itemize}
The first property follows from the fact that any element of
$\roots $ is lying in a basis of $\ndZ ^n$, and from
(\ref{eq-si}). The second one will be a corollary of Theorem
\ref{t-corr}.

\section{The Weyl--Brandt groupoid of a Nichols algebra of
diagonal type}
\label{sec-WBg}

Let $k$ be a field of characteristic zero, $G$ an abelian group,
and $V$ a \YD module over $kG$ of rank $n$ for some $n\in \ndN $.
Assume that $V$ is of diagonal type. More precisely, let
$\{x_1,x_2,\ldots ,x_n\}$ be a basis of $V$, $\{g_i\,|\,1\le i\le
n\}$ a subset of $G$, and $q_{ij}\in k\setminus \{0\}$ for all
$i,j\in \{1,2,\ldots ,n\}$, such that
\begin{align*}
\lcoa (x_i)=&g_i\ot x_i,& g_i\lact x_j=&q_{ij}x_j
\end{align*}
for $i,j\in \{1,2,\ldots ,n\}$. Then $V$ is a braided vector space
\cite[Def.\,5.4]{inp-Andr02} with braiding $\sigma \in \End (V\ot
V)$, where
\begin{align*}
\sigma (x_i\ot x_j)=q_{ij}x_j\ot x_i
\end{align*}
for all $i,j\in \{1,2,\ldots ,n\}$. The Nichols algebra $\cB (V)$
is of diagonal type \cite[Def.\,5.8]{inp-Andr02} and has a $\ndZ
^n$-grading defined by
\begin{align*}
\deg x_i:=\Ndbasis _i\quad \text{for all $i\in \{1,2,\ldots ,n\}$}
\end{align*}
where $E_0:=\{\Ndbasis _1,\ldots ,\Ndbasis _n\}$ is the standard
basis of the $\ndZ $-module $\ndZ ^n$. Let $\chi :\ndZ ^n\times
\ndZ ^n\to k\setminus \{0\}$ denote the bicharacter defined by
\begin{align}\label{eq-chi}
\chi (\Ndbasis _i,\Ndbasis _j):=q_{ij},
\end{align}
where $i,j\in \{1,2,\ldots ,n\}$. Kharchenko proved \cite[Theorem
2]{a-Khar99} that the algebra $\cB (V)$ has a (restricted) \PBW
basis consisting of iterated skew-commutators of the elements
$x_i$ of $V$ and satisfying the following property (see
\cite[Sect.~2]{a-Heck04d}):
\begin{itemize}
\item[(P)]
the height of a \PBW generator of $\ndZ ^n$-degree $d$ is finite
if and only if $2\le \mathrm{ord}\,\chi (d,d)<\infty $
($\mathrm{ord}$ means order with respect to multiplication), and
in this case it coincides with $\mathrm{ord}\,\chi (d,d)$.
\end{itemize}

In \cite[Sect.~3]{a-Heck04c} $\proots (\cB (V))$ was defined as
the set of degrees of the (restricted) \PBW generators counted
with multiplicities. The definition is independent of the choice
of a $\ndZ ^n$-graded \PBW basis satisfying property (P). Write
$\roots (\cB (V)):=\proots (\cB (V))\cup -\proots (\cB (V))$.
Since the elements of $\proots (\cB (V))$ are lying in $\ndN
_0^n\setminus \{0\}$, this is a disjoint union.

In \cite[Sect.~5]{a-Heck04c} the Weyl--Brandt groupoid $W(V)$
associated to $\cB (V)$ was defined as $W_{\chi ,E_0}$, where
$E_0$ is the standard basis of $\ndZ ^n$ and $\chi $ is defined by
(\ref{eq-chi}).

\begin{thm}\label{t-corr}
If $\chi :\ndZ ^n\times \ndZ ^n\to k\setminus \{0\}$ is a
bicharacter on $\ndZ ^n$ and $(\roots ,\chi ,E_0)$ is an
arithmetic root system then for the braided vector space $V$ of
diagonal type with $\dim V=n$ and with structure constants
$q_{ij}:=\chi (\Ndbasis _i,\Ndbasis _j)$ one has $\roots (\cB
(V))=\roots $. Conversely, if $V$ is a braided vector space of
diagonal type such that $\proots (\cB (V))$ is finite then
$(\roots (\cB (V)),\chi ,E_0)$ is an arithmetic root system, where
$\chi :\ndZ ^n\times \ndZ ^n\to k\setminus \{0\}$ is defined by
(\ref{eq-chi}). These correspondences are inverse to each other.
\end{thm}

First let us collect some consequences of this theorem.

\begin{folg}
For any arithmetic root system $(\roots ,\chi ,E)$, where $\chi $
has values in $k\setminus \{0\}$, the property $(AR2)$ is
fulfilled.
\end{folg}

\begin{folg}\label{f-multi}
For any braided vector space $V$ of diagonal type such that
$\roots (\cB (V))$ is finite, the multiplicities of the elements
of $\roots (\cB (V))$ are one, and if $\alpha \in \proots (\cB
(V))$ then $\lambda \alpha \notin \proots (\cB (V))$ for all
$\lambda \in \ndR \setminus \{1\}$.
\end{folg}

\begin{bew}[ of the theorem]
Let $\chi :\ndZ ^n\times \ndZ ^n\to k\setminus \{0\}$ be a
bicharacter on $\ndZ ^n$, $(\roots ,\chi ,E_0)$ an arithmetic root
system, and $V$ as in the theorem. By Definition \ref{d-ars} the
Weyl--Brandt groupoid $W_{\chi ,E_0}=W(V)$ is finite and for all
ordered bases $E$ of $\ndZ ^n$ with $(\id ,E)\in W(V)$ the maps
$s_{i,E}$ are well-defined for all $i\in \{1,2,\ldots ,n\}$.

By \cite[Prop.~1]{a-Heck04c} one has $\roots \subset \roots (\cB
(V))$. One has to show that

1. any $\alpha \in \roots $ has multiplicity 1 in $\roots (\cB
(V))$,

2. if $\alpha \in \roots (\cB (V))$, $\beta \in \roots $, $\lambda
\in \ndR $, and $\alpha =\lambda \beta $, then $\lambda ^2=1$, and

3. there is no $\alpha \in \roots (\cB (V))$ such that $\alpha
\notin \ndR \beta $ for all $\beta \in \roots $.

Assertions 1 and 2 follow from the fact that any $\alpha \in
\roots \subset \roots (\cB (V))$ is an element of some basis $E$
of $\ndZ ^n$ with $(\id ,E)\in W(V)$, and for the degrees of
generators of Nichols algebras these assertions are known.

Let $\Rnvec{n}$ be a vector in $(\ndR ^+)^n$, where $\ndR ^+=
\{r\in \ndR \,|\,r>0\}$, and let $\scalp{\cdot }{\cdot }$ denote
the standard scalar product on $\ndR ^n$. One clearly has the
equation
\begin{align}
\proots (\cB (V))&=\{\alpha \in \roots (\cB (V))\,|\,\scalp{\alpha
}{\Rnvec{n}}>0\}.
\end{align}
It is not difficult to check that there exists a continuous family
of vectors $\Rnvec{n}(t)$, where $t\in [0,1]$,
$\Rnvec{n}(0)=\Rnvec{n}$, and $\Rnvec{n}(1)=-\Rnvec{n}$, such that
$\Rnvec{n}(t)\not=0$ for all $t$ and the hyperplanes $H(t)$
containing 0 and orthogonal to $\Rnvec{n}(t)$ never contain more
than one line $\ndR\alpha $ with $\alpha \in \proots _E\subset
\ndZ ^n \subset \ndR ^n$. Additionally one can assume that there
exists a sequence $(\alpha _1,\alpha _2,\ldots ,\alpha _m)$ of
elements of $\proots _E$ (which is in general not unique) such
that

 1. for all $i\le m$ there exists a unique $t_i\in (0,1)$ such that
 $\alpha _i\in H(t_i)$,

 2. if $t\not=t_i$ for all $i\in \{1,2,\ldots ,m\}$ then $H(t)\cap
 \roots =\{\}$, and

 3. $i<j$ implies that $t_i<t_j$.

Take $\alpha \in \proots (\cB (V))$. Since $\scalp{\alpha
}{\Rnvec{n}(0)}>0$ and $\scalp{\alpha }{\Rnvec{n}(1)}<0$, there
exists $t\in (0,1)$ such that $\alpha \in H(t)$. In particular,
the first of the above requirements states that if additionally
$\alpha \in \proots _E$ then there exists a uniquely determined
$j\in \{1,2,\ldots ,m\}$ such that $\alpha = \alpha _j$. We will
show by induction that
\begin{itemize}
\item[($*$)]
for all $t\in [0,1]$ with $H(t)\cap \roots =\{\}$ there exists an
ordered basis $E(t)\subset \roots $ of $\ndZ ^n$ such that $(\id
,E(t))\in W(V)$, $\scalp{\Rnvec{n}(t)}{e_i(t)}>0$ for all
$e_i(t)\in E(t)$.
\end{itemize}
This gives the first assertion of the theorem. Indeed, relations
$\scalp{\Ndbasis _i}{\Rnvec{n}(1)}<0$ for all $i\in \{1,2,\ldots
,n\}$ imply that $(\id ,-E'_0)\in W(V)$, where $E'_0=E_0$ up to
permutation of its elements. Suppose that $\alpha \in \proots (\cB
(V))$. Then there exist ordered bases $E'$ and $E''$ of $\ndZ ^n$
and $i\in \{1,2,\ldots, n\}$, such that $(\id ,E'),(\id ,E'')\in
W(V)$, $E''=s_{i,E'}(E')$, and $\alpha $ is a nonnegative
respectively nonpositive integer linear combination of the
elements of $E'$ and $E''$, respectively. By the definition of
$s_{i,E'}$ at most one of the coefficients of $\alpha $ changes
with the base change $E'\to E''$. Hence $\alpha $ is a multiple of
an element in $E'$, that is $\alpha \in \ndR \beta $ for some
$\beta \in \roots $.

Assertion ($*$) clearly holds for $t<t_1$ with $E(t)=E_0$. Suppose
that it is valid for all $t\in [0,t_j)$ with some $j\le m$, and
set $E':=E(t)=\{e'_1,\ldots ,e'_n\}$ for some $t_{j-1}<t<t_j$.
First note that $\alpha _j\in E'$. Indeed, since
$\scalp{e'_i}{\alpha _j}>0$ for all $i\in \{1,2,\ldots ,n\}$, and
$\roots =\proots _{E'}\cup -\proots _{E'}$, one has at least
$\alpha _j\in \proots _{E'}$. Therefore from $\scalp{\alpha
_j}{\Rnvec{n}(t_j)}=0$ one concludes that also
$\scalp{e'_l}{\Rnvec{n}(t_j)}=0$ for some $l\in \{1,2,\ldots
,n\}$. Since $e'_l\in \roots $, uniqueness of $t_j$ gives $\alpha
_j=e'_l$.

By definition of $H(t)$, one has $\scalp{\Rnvec{n}(t)}{e'_l}<0$
for $t>t_j$, and if $\alpha \in \roots \setminus \ndR \alpha _j$,
$\scalp{\Rnvec{n}(t)}{\alpha }>0$ for $t_{j-1}<t<t_j$, then also
$\scalp{\Rnvec{n}(t)}{\alpha }>0$ for $t_j<t<t_{j+1}$ (where
$t_{m+1}:=1$). In particular,
$\scalp{\Rnvec{n}(t)}{s_{l,E'}(e'_i)}>0$ for all $i\in
\{1,2,\ldots ,n\}$ and $t_j<t<t_{j+1}$. Thus one can choose
$E(t)=s_{l,E'}(E')$ if $t_j<t<t_{j+1}$. This proves ($*$).

If $V$ is a braided vector space of diagonal type with $\dim V=n$
and $\proots (\cB (V))$ is finite, then $W(V)$ is finite and for
all ordered bases $E'$ of $\ndZ ^n$ with $(\id ,E')\in W(V)$ the
maps $s_{i,E'}$ are well-defined for all $i\in \{1,2,\ldots ,n\}$.
Hence $(\roots ',\chi ,E_0)$, where $\chi $ is defined by
(\ref{eq-chi}) and $\roots '=\bigcup \{E'\subset \ndZ ^n\,|\,(\id
,E')\in W_{\chi ,E_0}\}$, is an arithmetic root system. By the
first part of the theorem one gets $\roots (\cB (V))=\roots '$.
This proves the second part.

The last assertion follows from the fact that $\roots $ and
$\roots (\cB (V))$ are uniquely determined by $\chi $ and $V$,
respectively, and the latter are in one-to-one correspondence via
(\ref{eq-chi}).
\end{bew}

In \cite{a-Heck04d} the notion of Weyl equivalence of braided
vector spaces was introduced as follows. Let $V',V''$ be two rank
$n$ braided vector spaces of diagonal type. If with respect to
certain bases their structure constants $q'_{jl}$ and $q''_{jl}$,
satisfy the equations
\begin{align}
q'_{jj}=&q''_{jj},& q'_{jl}q'_{lj}=&q''_{jl}q''_{lj}
\end{align}
for all $j,l\in \{1,2,\ldots ,n\}$ then one says that \textit{$V'$
and $V''$ are twist equivalent}
\cite[Definition\,3.8]{inp-AndrSchn02}. Further, $V'$ and $V''$
are called \textit{Weyl equivalent}, if there exists a braided
vector space $V$ of diagonal type which is twist equivalent to
$V'$, and an ordered basis $E''$ of $\ndZ ^n$ such that $(\id
,E'')\in W(V'')$ and $V$ has degree $E''$ with respect to $V''$
\cite[Def.~1]{a-Heck04d}.

\begin{defin}
Let $E$ be an ordered basis of $\ndZ ^n$. Then two arithmetic root
systems $(\roots ',\chi ',E)$ and $(\roots '',\chi '',E)$ are
called \textit{Weyl equivalent} if there exist linear maps $\tau
,T\in \Aut (\ndZ ^n)$ such that $(T,E)\in W_{\chi '',E}$, $\tau
T(E)\subset T(E)$, $\roots ''=\tau T(\roots ')$, and $\chi
'(e,e)=\chi ''(\tau T(e),\tau T(e))$ for all $e\in \ndZ ^n$.
Further, they are called \textit{twist equivalent} if they are
Weyl equivalent with $T=\id $.
\end{defin}

\begin{bem} Let $V'$ and $V''$ be braided vector spaces of
diagonal type such that $\roots (\cB (V'))$ and $\roots (\cB
(V''))$ are finite. Then $V'$ and $V''$ are twist equivalent
respectively Weyl equivalent if and only if the corresponding
arithmetic root systems $(\roots (\cB (V')),\chi ',E_0)$ and
$(\roots (\cB (V'')),\chi '',E_0)$ have this property.
\end{bem}

Since for any given arithmetic root system one can easily
determine all arithmetic root systems which are Weyl equivalent to
it, it is sufficient to determine the Weyl equivalence classes.
This will be done for the rank 2 case in the next section.

\section{Arithmetic root systems of rank 2}
\label{sec-Ars2}

In \cite[Prop.~2]{a-Heck04d} Weyl equivalence of certain rank 2
braided vector spaces of diagonal type was considered. The main
result of the present paper is that the list given there, see
Figure 1, contains all rank 2 braided vector spaces $V$ of
diagonal type, such that $\roots (\cB (V))$ is finite.

\begin{figure}
 {\scriptsize
\begin{align*}
\begin{array}{r|l|l|l|c}
 & \text{conditions for $q_{ij}$; $q_{12}$ is always 1} & \text{free parameters} &
 \text{fixed parameters} & \text{tree}\\
\hline \hline
 1 & q_{21}=1 & & q_{11},q_{22}\in k\setminus \{0\} & T1\\
\hline
 2 & q_{21}=q_{11}^{-1},\ q_{22}=q_{11} & & q_{11}\in k\setminus \{0,1\} & T2\\
\hline
 3 & q_{11}=q,\ q_{21}=q^{-1},\ q_{22}=-1,\text{ or} & &
 q\in k\setminus \{-1,0,1\} & T2\\
 & q_{11}=-1,\ q_{21}=q,\ q_{22}=-1 & & & T2\\
\hline
 4 & q_{11}=q,\ q_{21}=q^{-2},\ q_{22}=q^2 & & q\in k\setminus \{-1,0,1\} & T3\\
\hline
 5 & q_{11}=q,\ q_{21}=q^{-2},\ q_{22}=-1 & q\in \{q_0,-q_0^{-1}\} &
 q_0\in k\setminus \{-1,0,1\} & T3\\
\hline
 6 & q_{11}=\zeta ,\ q_{21}=q^{-1},\ q_{22}=q & q\in \{q_0,\zeta q_0^{-1}\} &
 \zeta \in R_3 & T3\\
 & & & q_0\in k\setminus \{0,1,\zeta ,\zeta ^2\} &\\
\hline
 7 & q_{11}=\zeta ,\ q_{21}=-\zeta ,\ q_{22}=-1 & \zeta \in R_3 & & T3\\
\hline
 8 & q_{11}=\zeta ^4,\ q_{21}=\zeta ^{-3},\ q_{22}=-\zeta ^2,\text{ or}&
 \zeta \in \{\zeta _0,-\zeta _0^{-1}\}&
 \zeta _0\in R_{12} & T4\\
 & q_{11}=\zeta ^4,\ q_{21}=\zeta ^{-1},\ q_{22}=-1,\text{ or}& & & T5\\
 & q_{11}=\zeta ^{-3},\ q_{21}=\zeta ,\ q_{22}=-1 & & & T7\\
\hline
 9 & q_{11}=-\zeta ^2,\ q_{21}=\zeta ,\ q_{22}=-\zeta ^2,\text{ or}&
& \zeta \in R_{12} & T4\\
 & q_{11}=-\zeta ^2,\ q_{21}=\zeta ^3,\ q_{22}=-1,\text{ or}& & & T5\\
 & q_{11}=-\zeta ^{-1},\ q_{21}=\zeta ^{-3},\ q_{22}=-1 & & & T7\\
\hline
 10 & q_{11}=\zeta ,\ q_{21}=\zeta ^{-2},\ q_{22}=-\zeta ^3,\text{ or}&
& \zeta \in R_{18} & T6\\
 & q_{11}=-\zeta ^2,\ q_{21}=-\zeta ,\ q_{22}=-1,\text{ or} & & & T14\\
 & q_{11}=-\zeta ^3,\ q_{21}=-\zeta ^{-1},\ q_{22}=-1 & & & T9\\
\hline
 11 & q_{11}=q,\ q_{21}=q^{-3},\ q_{22}=q^3 & & q\in k\setminus \{-1,0,1\} & T8\\
 & & & q\notin R_3 & \\
\hline
 12 & q_{11}=\zeta ^2,\ q_{21}=\zeta ,\ q_{22}=\zeta ^{-1},\text{ or}&
& \zeta \in R_8 & T8\\
 & q_{11}=\zeta ^2,\ q_{21}=-\zeta ^{-1},\ q_{22}=-1,\text{ or} & & & T8\\
 & q_{11}=\zeta ,\ q_{21}=-\zeta ,\ q_{22}=-1 & & & T8\\
\hline
 13 & q_{11}=\zeta ^6,\ q_{21}=-\zeta ^{-1},\ q_{22}=\zeta ^8,\text{ or} & &
\zeta \in R_{24} & T10\\
 & q_{11}=\zeta ^6,\ q_{21}=\zeta ,\ q_{22}=\zeta ^{-1},\text{ or}& & & T13\\
 & q_{11}=\zeta ^8,\ q_{21}=\zeta ^5,\ q_{22}=-1,\text{ or} & & & T17\\
 & q_{11}=\zeta ,\ q_{21}=\zeta ^{-5},\ q_{22}=-1 & & & T21\\
\hline
 14 & q_{11}=\zeta ,\ q_{21}=\zeta ^{-3},\ q_{22}=-1,\text{ or} &
\zeta \in \{\zeta _0,\zeta _0^{11}\}& \zeta _0\in R_5\cup R_{20} & T11\\
 & q_{11}=-\zeta ^{-2},\ q_{21}=\zeta ^3,\ q_{22}=-1 & & & T16\\
\hline
 15 & q_{11}=\zeta ,\ q_{21}=\zeta ^{-3},\ q_{22}=-\zeta ^5,\text{ or} & &
\zeta \in R_{30} & T12\\
 & q_{11}=-\zeta ^3,\ q_{21}=-\zeta ^4,\ q_{22}=-\zeta ^{-4},\text{ or}& & & T15\\
 & q_{11}=-\zeta ^5,\ q_{21}=-\zeta ^{-2},\ q_{22}=-1,\text{ or} & & & T18\\
 & q_{11}=-\zeta ^3,\ q_{21}=-\zeta ^2,\ q_{22}=-1 & & & T20\\
\hline
 16 & q_{11}=\zeta ,\ q_{21}=\zeta ^{-3},\ q_{22}=-1,\text{ or} & &
\zeta \in R_{14} & T19\\
 & q_{11}=-\zeta ^{-2},\ q_{21}=\zeta ^3,\ q_{22}=-1 & & & T22\\
\hline
\end{array}
\end{align*}
} \caption{Weyl equivalence for rank 2 \YD modules}
\end{figure}

\begin{thm}\label{t-class2}
Let $V$ be a braided vector space of diagonal type with $\dim
V=2$. Then $\roots (\cB (V))$ is finite if and only if there
exists a rank 2 braided vector space of diagonal type which is
twist equivalent to $V$ and has structure constants $q_{ij}$ which
appear in Figure 1.
\end{thm}

For the classification of arithmetic root systems of rank~2
several technical lemmata will be needed. The first one is a
well-known classical result, which can be easily proved with help
of eigenvalue considerations.

\begin{lemma}\label{l-finord}
Let $A\in \SLZ{2}$. Then $\mathrm{ord}\,A$ is finite if and only
if either $A=\id $ or $A=-\id $ or $\det A=1$, $\tr A\in
\{-1,0,1\}$.
\end{lemma}

\begin{lemma}\label{l-subSLZ}
For a given $M\in \ndN $ let $S$ be a subsemigroup of $\SLZ{2}$
such that its generators are of the form $\begin{pmatrix}a & -b\\
c & -d\end{pmatrix}$ with $0<d<Mb<a$. Then all elements of $S$ are
of this form. In particular, $\id \notin S$.
\end{lemma}

\begin{bew}
Assume that $A_1:=\begin{pmatrix}a_1 & -b_1\\
c_1 & -d_1\end{pmatrix}\in S$ and $A_2:=\begin{pmatrix}a_2 & -b_2\\
c_2 & -d_2\end{pmatrix}\in S$ such that $0<d_1<Mb_1<a_1$ and
$0<d_2<Mb_2<a_2$. Since $\det A_1=\det A_2=1$, one also has
$c_1=(a_1d_1+1)/b_1$ and $c_2=(a_2d_2+1)/b_2$. Then
$A_1A_2=\begin{pmatrix}a_1a_2-b_1c_2 & -(a_1b_2-b_1d_2) \\
c_1a_2-d_1c_2 & -(c_1b_2-d_1d_2)\end{pmatrix}$. We obtain the
following inequalities.
\begin{align}
\label{ineq1}
 a_1b_2-b_1d_2=&b_2(a_1-Mb_1)+b_1(Mb_2-d_2)>0,\\
\label{ineq2}
 c_1b_2-d_1d_2=&\frac{a_1d_1+1}{b_1}b_2-d_1d_2=
\frac{d_1(a_1b_2-b_1d_2)+b_2}{b_1}>0 \text{ by (\ref{ineq1}),}
\end{align}
\begin{align}
  &M(a_1b_2-b_1d_2)-(c_1b_2-d_1d_2)=(Mb_1-d_1)(Mb_2-d_2) \notag \\
  &\qquad +b_2\left(Ma_1-M^2b_1+Md_1
  -\frac{a_1d_1+1}{b_1}\right) \label{ineq3}\\
  &\qquad
  =(Mb_1-d_1)(Mb_2-d_2)+\frac{b_2}{b_1}((Mb_1-d_1)(a_1-Mb_1)-1)>0, \notag
\\
  &a_1a_2-b_1c_2-M(a_1b_2-b_1d_2)=(a_1-Mb_1)(a_2-Mb_2)\notag \\
  &\qquad +b_1\left(Ma_2-M^2b_2-\frac{a_2d_2+1}{b_2}+Md_2\right)
  \label{ineq4}\\
  &\qquad=(a_1-Mb_1)(a_2-Mb_2)+\frac{b_1}{b_2}((Mb_2-d_2)(a_2-Mb_2)-1)>0.
  \notag
\end{align}
The inequalities (\ref{ineq2})--(\ref{ineq4}) give the
assertion.
\end{bew}

For any ordered subset $E=\{e_1,e_2\}$ of $\ndZ ^2$ consisting of
two elements let $\tau (E)$ denote the transposed set
$\{e_2,e_1\}$. A useful consequence of Definition \ref{d-ars} is
the following.

\begin{lemma}\label{l-seq}
Let $\chi $ be a bicharacter on $\ndZ ^2$ with values in
$k\setminus \{0\}$. For all $i\in \ndN _0$ set $T'_i=s_{1,E'_i}$
and $T''_i=s_{1,E''_i}$, where $E'_{i+1}=\tau (T'_i(E'_i))$,
$E''_{i+1}=\tau (T''_i(E''_i))$, and $E'_0=E_0$, $E''_0=\tau
(E_0)$. The triple $(\roots ,\chi ,E_0)$, where $\roots =\bigcup
\{E\subset \ndZ ^2\,|\,(\id,E)\in W_{\chi ,E_0}\}$, is an
arithmetic root system if and only if the sequences
$(T'_iT'_{i-1}\cdots T'_0)_{i\in \ndN _0}$ and
$(T''_iT''_{i-1}\cdots T''_0)_{i\in \ndN _0}$ of elements of $\Aut
(\ndZ ^2)$ are well-defined and periodic. In particular, they
contain the identity.
\end{lemma}

\begin{bew}
The if part is clear, since the equations
$s_{i,s_{i,E}(E)}s_{i,E}=\id $ hold for all $E$ with $(\id ,E)\in
W_{\chi ,E_0}$. On the other hand, if $W(V)$ is finite, then the
set $\{E'_i\,|\,i\in \ndN _0\}$ is finite. Since $T'_i$ is
invertible and it depends only on $E'_i$, one obtains that the
sequence $(T'_i)_{i\in \ndN _0}$ is periodic. Let $j\in \ndN $ be
the smallest number such that $E'_j=E_0$. Then $T'_i=T'_{i+j}$ for
all $i\in \ndN _0$, and hence $T'_{ij-1}T'_{ij-2}\cdots T'_0 =
(T'_{j-1}\cdots T'_0)^i$ for all $i\in \ndN $. Again, the set
$\{(T'_{j-1}\cdots T'_0)^i(E_0)\,|\,i\in \ndN _0\}$ is finite, and
hence there exist $i_1<i_2\in \ndN _0$ such that $(T'_{j-1}\cdots
T'_0)^{i_1}(E_0) = (T'_{j-1}\cdots T'_0)^{i_2}(E_0)$, that is
$(T'_{j-1}\cdots T'_0)^{i_2-i_1}(E_0) = E_0$. The proof for the
second sequence is similar.
\end{bew}

\begin{lemma}\label{l-no1}
Let $V$ be a braided vector space of diagonal type with $\dim
V=2$, and let $(q_{ij})_{i,j\in \{1,2\}}$ be the matrix of its
structure constants with respect to a certain basis. If $\roots
(\cB (V))$ is finite then either
\begin{gather*}
 (q_{12}q_{21}-1)(q_{11}q_{12}q_{21}-1)(q_{12}q_{21}q_{22}-1)
 (q_{11}+1)(q_{22}+1)=0\qquad \text{or}\\
 q_{11}q_{12}^2q_{21}^2q_{22}=-1,\
 (q_{11}^2+q_{11}+1)(q_{11}^2+1)(q_{11}^2q_{12}q_{21}-1)(q_{11}^3q_{12}q_{21}-1)=0.
\end{gather*}
\end{lemma}

\begin{bew}
Suppose that $\roots (\cB (V))$ is finite. By Corollary
\ref{f-multi}, $\proots (\cB (V))$ doesn't contain $2(\Ndbasis
_1+\Ndbasis _2)$ and $2(2\Ndbasis _1+\Ndbasis _2)$. Set
$z_1:=x_1x_2-q_{12}x_2x_1$, $z_2:=x_1z_1-q_{11}q_{12}z_1x_1$, and
$z_3:=x_1z_2-q_{11}^2q_{12}z_2x_1$, and assume that
$(q_{12}q_{21}-1)(q_{11}q_{12}q_{21}-1)(q_{12}q_{21}q_{22}-1)
(q_{11}+1)(q_{22}+1)\not=0$. Then $z_1\not=0$ and $z_2\not=0$ by
\cite[Lemma~3.7]{inp-AndrSchn02} or \cite[Sect.~4.1]{a-Heck04b}.
Consider first the case when $\chi (\deg z_1,\deg z_1)=-1$. By
\cite[Sect.~4.2]{a-Heck04b} $z_1^{2m}\not=0$ for all $m\in \ndN $
and hence $\mathrm{ord}\,\chi (\deg z_1,\deg z_1)=2$ implies that
$2(\Ndbasis _1+\Ndbasis _2)\in \proots (\cB (V))$, which is a
contradiction. Similarly, if $\chi (\deg z_2,\deg z_2)=-1$ then
either $z_3=0$ (and hence $q_{11}\in R_3$, since the Equation
$q_{22}=-1$, which follows from $\chi (\deg z_2,\deg z_2)=-1$ and
$q_{11}^2q_{12}q_{21}=1$, was excluded), which corresponds to a
conclusion of the lemma, or $z_2^m\not=0$ for all $m\in \ndN $ by
\cite[Sect.~4.2]{a-Heck04b}, which is a contradiction to
$\mathrm{ord}\,\chi (\deg z_2,\deg z_2)=2$ and $2(2\Ndbasis
_1+\Ndbasis _2)\notin \proots (\cB (V))$.

Assume now that $\chi (\deg z_i,\deg z_i)\not=-1$ for $i\in
\{1,2\}$. Since $2(\Ndbasis _1+\Ndbasis _2)\notin \proots (\cB
(V))$, the skew-commutator $z_2x_1-\chi (\deg z_2,\deg x_1)x_1z_2$
must not be a \PBW generator of $\cB (V)$. By
\cite[Sect.~4.4]{a-Heck04b} this gives that $w_1=0$ (in the
notation there) and hence $q_{11}q_{12}^2q_{21}^2q_{22}=-1$. If
$z_3=0$ then the proof is complete. Otherwise, since $2(2\Ndbasis
_1+\Ndbasis _2)\notin \proots (\cB (V))$, the skew-commutator
$z_3x_1-\chi (\deg z_2,\deg x_1)x_1z_2$ must not be a \PBW
generator of $\cB (V)$. Again by \cite[Sect.~4.4]{a-Heck04b} this
gives that $w_2=0$ and hence
\begin{align*}
-q_{11}^{-1}q_{21}^{-3}\frac{(1+q_{11}^{-2}) (1-q_{12}q_{21})^2
(1-q_{11}^3q_{12}q_{21})}{1-q_{11}^{-1}}z_3=0.
\end{align*}
This proves the lemma.
\end{bew}

Assign to $V$ two sequences $(m'_i)_{i\in \ndN _0}$ and
$(m''_i)_{i\in \ndN _0}$, such that $m'_i=m_{12}$ and
$m''_i=m_{12}$ as defined in (\ref{eq-mij}) with $E=E'_i$ and
$E=E''_i$, respectively (in the notation of Lemma \ref{l-seq}).
Further, for any $l\in \ndN _0$ set $\scq{'}{ij}{l}:=\chi
(e'_i,e'_j)$ and $\scq{''}{ij}{l}:=\chi (e''_i,e''_j)$, where
$E'_l=\{e'_1,e'_2\}$ and $E''_l=\{e''_1,e''_2\}$, and define
\begin{align}\label{eq-pprime}
 \scp{'}{l}:=&
 \begin{cases}
  1 & \text{if $(\scq{'}{11}{l})^{m'_l}\scq{'}{12}{l}\scq{'}{21}{l}=1$,}\\
  (\scq{'}{11}{l})^{-1}\scq{'}{12}{l}\scq{'}{21}{l} &
  \text{otherwise,}
 \end{cases}\\ \label{eq-ppprime}
 \scp{''}{l}:=&
 \begin{cases}
  1 & \text{if $(\scq{''}{11}{l})^{m''_l}\scq{''}{12}{l}\scq{''}{21}{l}=1$,}\\
  (\scq{''}{11}{l})^{-1}\scq{''}{12}{l}\scq{''}{21}{l} &
  \text{otherwise.}
 \end{cases}
\end{align}
By Equations (3) in \cite{a-Heck04d} one obtains for all $l\in
\ndN _0$ the relations
\begin{align}\label{eq-qprime}
\scq{'}{11}{l+1}=&\scp{'}{l}{}^{m'_l}\scq{'}{22}{l},&
\scq{'}{12}{l+1}\scq{'}{21}{l+1}=&\scp{'}{l}{}^{-2}
\scq{'}{12}{l}\scq{'}{21}{l},& \scq{'}{22}{l+1}=&\scq{'}{11}{l},\\
\label{eq-qpprime}
\scq{''}{11}{l+1}=&\scp{''}{l}{}^{m''_l}\scq{''}{22}{l},&
\scq{''}{12}{l+1}\scq{''}{21}{l+1}=&\scp{''}{l}{}^{-2}
\scq{''}{12}{l}\scq{''}{21}{l},&
\scq{''}{22}{l+1}=&\scq{''}{11}{l}.
\end{align}

\begin{bew}[ of Theorem \ref{t-class2}]
Let $(q_{ij})_{i,j\in \{1,2\}}$ be the matrix of structure
constants of $V$ (with respect to a given basis). Since it is
sufficient to determine Weyl equivalence classes, one can assume
that $q_{12}=1$. If $V$ appears in Figure 1 then $\roots (\cB
(V))$ is finite \cite[Sect.~3]{a-Heck04d}. It remains to show that
for any Weyl equivalence class at least one representant $V$, such
that $\roots (\cB (V))$ is finite, appears in Figure 1. We use
Lemma \ref{l-seq} to obtain restrictions on the structure
constants $q_{ij}$.

For all $i\in \ndN _0$ let $\tilde{T}'_i$ and $\tilde{T}''_i$
denote the matrix of $T'_i$ and $T''_i$, respectively, with
respect to the basis $E'_i$ and $E''_i$, respectively, and set
$\tilde{\tau }=\begin{pmatrix} 0 & 1 \\ 1 & 0 \end{pmatrix}$. Then
the matrix of $T'_iT'_{i-1}\cdots T'_0$ with respect to the basis
$E_0$ is $\tilde{T}'_0\tilde{\tau }\tilde{T}'_1\tilde{\tau }\cdots
\tilde{T}'_i$ if $i$ is even and $\tilde{T}'_0\tilde{\tau
}\tilde{T}'_1\tilde{\tau }\cdots \tilde{T}'_i\tilde{\tau }$ if $i$
is odd. Similarly, the matrix of $T''_iT''_{i-1}\cdots T''_0$ with
respect to the basis $E_0$ is $\tilde{\tau } \tilde{T}''_0
\tilde{\tau } \tilde{T}''_1 \tilde{\tau } \cdots \tilde{T}''_i
\tilde{\tau }$ if $i$ is even and $\tilde{\tau } \tilde{T}''_0
\tilde{\tau } \tilde{T}''_1 \tilde{\tau } \cdots \tilde{T}''_i$ if
$i$ is odd.

Using Weyl equivalence, one can assume that $m'_0$ is minimal
among all $m'_i$ and $m''_i$, where $i\in \ndN _0$.
If $m'_0=0$ then $q_{21}=1$ by (\ref{eq-mij}),
and then $V$ appears in the first line of Figure 1. We have to
consider now several cases.

Step 1. \textit{Assume that $m'_0>1$.} Then $m'_i>1$ for all $i\in
\ndN _0$. By Lemma \ref{l-seq} one of the maps
$T'_{2i+1}T'_{2i}\cdots T'_0$, where $i\in \ndN _0$, has to be the
identity. For all $j\in \ndN _0$ one has
\begin{align}
\tilde{T}_j\tilde{\tau }=&
\begin{pmatrix} -1 & m'_j \\ 0 & 1\end{pmatrix}
\begin{pmatrix} 0 & 1 \\ 1 & 0\end{pmatrix}
=\begin{pmatrix} m'_j & -1 \\ 1 & 0\end{pmatrix},
\end{align}
and hence $\tilde{T}_{2i} \tilde{\tau } \tilde{T}_{2i+1}
\tilde{\tau }=\begin{pmatrix} m'_{2i}m'_{2i+1}-1 & -m'_{2i} \\
m'_{2i+1} & -1\end{pmatrix}$ for all $i\in \ndN _0$. Since
$m'_j\ge 2$ for all $j\in \ndN _0$, these matrices satisfy the
conditions of Lemma \ref{l-subSLZ} with $M=1$. Then Lemma
\ref{l-subSLZ} gives that $T'_{2i+1}T'_{2i}\cdots T'_0$ is never
the identity, which is a contradiction.

Step 2. \textit{If $m'_0=1$ and $m''_0=1$} then one has either
$q_{11}q_{21}=q_{21}q_{22}=1$ or $q_{11}q_{21}=1$, $q_{22}=-1$, or
$q_{11}=-1$, $q_{21}q_{22}=1$ (which is twist equivalent to the
previous case), or $q_{11}=q_{22}=-1$. In all cases $V$ appears in
Figure 1.

Step 3. \textit{Suppose that $m'_{2i}=m''_{2i+1}=1$ and
$m'_{2i+1}>1$, $m''_{2i}>1$ for all $i\in \ndN _0$.} Again by Weyl
equivalence one can assume that $m'_1$ is minimal among all
$m'_{2i+1}$ and $m''_{2i}$ for all $i\in \ndN _0$. Recall that one
has $\tilde{T}_{2i} \tilde{\tau } \tilde{T}_{2i+1} \tilde{\tau }
=\begin{pmatrix} m'_{2i+1}-1 & -1
\\ m'_{2i+1} & -1\end{pmatrix}$, and hence if $m'_1>3$ (that is
$m'_{2i+1}\ge 4$ for all $i\in \ndN _0$) then these matrices
satisfy the conditions of Lemma \ref{l-subSLZ} with $M=2$. This
gives again a contradiction. Thus one has $m'_1\in \{2,3\}$, and
since $m'_0=1$, also $q_{11}q_{21}=1$ (see Step 3a) or
$q_{11}=-1$, $q_{21}^2\not=1$ (Step 3b) holds by (\ref{eq-mij}).

Step 3a. \textit{Assume that additionally $q_{11}q_{21}=1$ holds.}
Then (\ref{eq-qprime}) gives that $\scp{'}{0}=1$,
$\scq{'}{11}{1}=q_{22}$,
$\scq{'}{12}{1}\scq{'}{21}{1}=q_{11}^{-1}$,
$\scq{'}{22}{1}=q_{11}$. If $m'_1=2$ then (\ref{eq-mij}) gives
$q_{11}=q_{22}^2$ or $q_{22}\in R_3$. These cases appear in lines
4 and 6, respectively, of Figure 1, if one replaces $E_0$ by $\tau
(E_0)$, and in previous lines for certain special values of
$q_{11}$. If $m'_1=3$ then (\ref{eq-mij}) gives that either
$q_{11}=q_{22}^3$ or $q_{22}^2=-1$, $q_{11}^4\not=1$. In the first
case $V$ appears in Figure 1. In the second, from
(\ref{eq-qprime}) one gets $\scp{'}{1}=q_{11}^{-1}q_{22}^{-1}$,
$\scq{'}{11}{2}=q_{11}^{-2}q_{22}$,
$\scq{'}{12}{2}\scq{'}{21}{2}=-q_{11}$, $\scq{'}{22}{2}=q_{22}$.
Since $m'_2=1$ and $q_{11}^4\not=1$, (\ref{eq-mij}) yields
$\scq{'}{11}{2}=-1$, that is $q_{22}=-q_{11}^2$. This implies that
$q_{11}^4=-1$, and hence $V$ appears in the $12^\mathrm{th}$ line
of Figure 1.

Step 3b. \textit{Consider the case $q_{11}=-1$, $q_{21}^2\not=1$.}
Then
\begin{align*}
m'_0=&1,& \scp{'}{0}=&-q_{21},& \scq{'}{11}{1}=&-q_{21}q_{22},&
\scq{'}{12}{1}\scq{'}{21}{1}=&q_{21}^{-1},& \scq{'}{22}{1}=&-1.
\end{align*}
If $m'_1=2$ then either $q_{21}q_{22}^2=1$ or $-q_{21}q_{22}\in
R_3$, $q_{21}^3\not=1$. The first case appears in the
$5^\mathrm{th}$ line of Figure 1. In the second one obtains
\begin{align*}
\scp{'}{1}=&-q_{21}^{-2}q_{22}^{-1},&
\scq{'}{11}{2}=&q_{21}^{-1}q_{22},&
\scq{'}{12}{2}\scq{'}{21}{2}=&-q_{22}^{-1},&
\scq{'}{22}{2}=&-q_{21}q_{22}.
\end{align*}
{}From $m'_2=1$ and from $q_{21}^2\not=1$ one gets
$q_{21}=-q_{22}$. Since $-q_{21}q_{22}\in R_3$ and $q_{21}\notin
R_3$, this yields $q_{22}\in R_3$. This example appears in line 7
of Figure 1.

 If $m'_1=3$ then one has again two possibilities:
either $-q_{21}^2q_{22}^3=1$ or $q_{21}^2q_{22}^2=-1$,
$q_{21}^4\not=1$ (see Step 3b3 for the latter). In the first case
\begin{align*}
m'_1=&3,& \scp{'}{1}=&1,& \scq{'}{11}{2}=&-1,&
\scq{'}{12}{2}\scq{'}{21}{2}=&q_{21}^{-1},&
\scq{'}{22}{2}=&-q_{21}q_{22},\\
m'_2=&1,& \scp{'}{2}=&-q_{21}^{-1},& \scq{'}{11}{3}=&q_{22},&
\scq{'}{12}{3}\scq{'}{21}{3}=&q_{21},& \scq{'}{22}{3}=&-1.
\end{align*}
By (\ref{eq-mij}) one has to have $q_{22}^{m'_3}q_{21}=1$ (Step
3b1) or $q_{22}\in R_{m'_3+1}$, $q_{21}^{m'_3+1}\not=1$, where
$m'_3\ge 3$ (Step 3b2).

Step 3b1. \textit{The setting $-q_{21}^2q_{22}^3=1$,
$q_{22}^{m'_3}q_{21}=1$}. One has
\begin{align*}
\scp{'}{3}=&1,& \scq{'}{11}{4}=&-1,&
\scq{'}{12}{4}\scq{'}{21}{4}=&q_{21},& \scq{'}{22}{4}=&q_{22}.
\end{align*}
Thus the sequence $(m'_i)_{i\in \ndN _0}$ has period 4, and
\begin{align*}
\prod _{i=0}^3\tilde{T}_i\tilde{\tau }=
\begin{pmatrix} 1 & -1 \\ 1 & 0 \end{pmatrix}
\begin{pmatrix} 3 & -1 \\ 1 & 0 \end{pmatrix}
\begin{pmatrix} 1 & -1 \\ 1 & 0 \end{pmatrix}
\begin{pmatrix} m'_3 & -1 \\ 1 & 0 \end{pmatrix}=
\begin{pmatrix} m'_3-2 & -1 \\ 2m'_3-3 & -2 \end{pmatrix}.
\end{align*}
By Lemma \ref{l-seq} this matrix has to have finite order. By
Lemma \ref{l-finord} the latter happens if and only if $m'_3\in
\{3,4,5\}$. If $m'_3=3$ then equations $q_{22}^{m'_3}q_{21}=1$ and
$-q_{21}^2q_{22}^3=1$ imply that $q_{21}=q_{22}^{-3}$ and
$-q_{22}^3=1$ which is a contradiction to $q_{21}^2\not=1$. If
$m'_3=4$ then one gets $q_{21}=q_{22}^{-4}$ and $q_{22}^5=-1$.
Since $q_{21}^2\not=1$ this yields $q_{22}\in R_{10}$. This is the
example in line 14 of Figure 1 (with $\zeta _0\in R_5$). If
$m'_3=5$ then $q_{21}=q_{22}^{-5}$ and $-q_{22}^7=1$. Again
$q_{21}^2\not=1$ implies that $q_{22}\in R_{14}$. This example
appears in line 16 of Figure 1.

Step 3b2. \textit{The setting $-q_{21}^2q_{22}^3=1$, $q_{22}\in
R_{m'_3+1}$, $q_{21}^2\not=1$, $q_{21}^{m'_3+1}\not=1$, where
$m'_3\ge 3$.} One obtains
\begin{align*}
\scp{'}{3}=&q_{21}q_{22}^{-1},&
\scq{'}{11}{4}=&-q_{21}^{m'_3}q_{22},&
\scq{'}{12}{4}\scq{'}{21}{4}=&q_{21}^{-1}q_{22}^2,&
\scq{'}{22}{4}=&q_{22}.
\end{align*}
Since $m'_4=1$, one has either $q_{21}^{m'_3-3}=1$ (Step 3b2a) or
$q_{21}^{m'_3}q_{22}=1$, $q_{21}^{-1}q_{22}^2\not=-1$ (Step 3b2b).

Step 3b2a. \textit{As in Step 3b2, but additionally one has
$q_{21}^{m'_3-3}=1$}. Then
\begin{align*}
m'_4=&1,& \scp{'}{4}=&1,& \scq{'}{11}{5}=&q_{22},&
\scq{'}{12}{5}\scq{'}{21}{5}=&q_{21}^{-1}q_{22}^2,&
\scq{'}{22}{5}=&-q_{21}^3q_{22},\\
m'_5=&m'_3,& \scp{'}{5}=&q_{21}^{-1}q_{22},& \scq{'}{11}{6}=&-1,&
\scq{'}{12}{6}\scq{'}{21}{6}=&q_{21},& \scq{'}{22}{6}=&q_{22}.
\end{align*}
Thus the sequence $(m'_i)_{i\in \ndN _0}$ has period 6. One has
\begin{align*}
\prod _{i=0}^5\tilde{T}_i\tilde{\tau }= \prod _{i=0}^5
\begin{pmatrix} m'_i & -1 \\ 1 & 0 \end{pmatrix}=
\begin{pmatrix} m'_3{}^2-4m'_3+2 & 3-m'_3 \\
 2m'_3{}^2-7m'_3+3 & 5-2m'_3\end{pmatrix}
\end{align*}
By Lemma \ref{l-seq} this matrix has to have finite order. By
Lemma \ref{l-finord} if the latter happens then $-2\le
(m'_3-3)^2-2\le 2$, that is $(3\le )m'_3\le 5$. If $m'_3=3$ then
the above matrix has trace $-2$, but it is not equal to $-\id $,
and hence it doesn't have finite order. If $m'_3\in \{4,5\}$ then
the relations $q_{21}^{m'_3-3}=1$ and $q_{21}^2\not=1$ contradict
to each other.

Step 3b2b. \textit{As in Step 3b2, but additionally one has
$q_{21}^{m'_3}q_{22}=1$, $q_{21}^{-1}q_{22}^2\not=-1$.} Inserting
the equations $q_{21}^{m'_3}=q_{22}^{-1}$ and
$q_{22}^{m'_3}=q_{22}^{-1}$ into $(-q_{21}^2q_{22}^3)^{2m'_3}=1$
one obtains $q_{22}^{10}=1$. Since $q_{22}\in R_{m'_3+1}$ and
$m'_3\ge 3$, this means that $m'_3\in \{4,9\}$. If $m'_3=4$ then
$q_{22}\in R_5$, $q_{22}=q_{21}^{-4}$, and hence
$-q_{21}^2q_{22}^3=1$ implies that $q_{21}^{10}=-1$ and
$q_{21}^2\not=-1$. Thus $q_{21}\in R_{20}$. This example appears
in line 14 of Figure 1 with $\zeta _0\in R_{20}$. Finally, if
$m'_3=9$ then equations $q_{22}=q_{21}^{-9}$ and
$-q_{21}^2q_{22}^3=1$ imply that $q_{21}^{25}=-1$. Since
$q_{22}^5=-1$, one also has $q_{21}^{45}=-1$, that is
$q_{21}^5=-1$ and hence $q_{22}=q_{21}$. This is a contradiction
to $q_{21}^{m'_3+1}\not=1$.

Step 3b3. \textit{The setting $q_{11}=-1$, $q_{21}^2q_{22}^2=-1$,
$q_{21}^4\not=1$.} One has
\begin{align*}
m'_1=&3,& \scp{'}{1}=&q_{22},& \scq{'}{11}{2}=&-q_{22}^3,&
\scq{'}{12}{2}\scq{'}{21}{2}=&-q_{21},&
\scq{'}{22}{2}=&-q_{21}q_{22}.
\end{align*}
Since $m'_2=1$, Equation (\ref{eq-mij}) gives that either
$q_{21}q_{22}^3=1$ or $q_{22}^3=1$. In the first case
$q_{21}^2q_{22}^2=-1$ gives that $q_{21}=-q_{22}$ and $q_{22}\in
R_8$. This example appears in line 12 of Figure 1. In the second
case again by Equation $q_{21}^2q_{22}^2=-1$ one concludes that
$q_{22}=-q_{21}^2$ and hence $q_{21}^6=-1$. Since $q_{21}^4\not=1$
this yields $q_{21}\in R_{12}$, which corresponds exactly to an
example in line 8 of Figure 1. (But this is an accident. One can
show that the conditions of Step 3 are not fulfilled. We will meet
this example another time.)

Step 4. If the setting is \textit{different from steps 1--3} then
there exists an $i\in \ndN _0$ such that $m'_i=1$ and
$m'_{i+1}>1$, $m'_{i+2}>1$. By Weyl equivalence one can assume
that $m''_0=1$ and $m'_0>1$, $m'_1>1$. Equations $m''_0=1$ and
(\ref{eq-mij}) give that either $q_{21}q_{22}=1$ (Step 4a) or
$q_{22}=-1$, $q_{21}^2\not=1$ (Step 4b).

Step 4a. \textit{The setting $q_{21}q_{22}=1$, $m'_0>1$,
$m'_1>1$.} By (\ref{eq-mij}) one has either
$q_{11}^{m'_0}q_{21}=1$ or $q_{11}\in R_{m'_0+1}$,
$q_{21}^{m'_0+1}\not=1$. In the first case $p'_0=1$ and $m'_1=1$
which is a contradiction. In the second one has
\begin{align}\label{eq-4a-1}
\scp{'}{0}=&q_{11}^{-1}q_{21},&
\scq{'}{11}{1}=&q_{11}q_{21}^{m'_0-1},&
\scq{'}{12}{1}\scq{'}{21}{1}=&q_{11}^2q_{21}^{-1},&
\scq{'}{22}{1}=&q_{11}.
\end{align}
Lemma \ref{l-no1} and the inequalities $m'_0>1$ and $m'_1>1$ imply
that $q_{11}^6q_{21}^{m'_0-3}=-1$. By (\ref{eq-mij}) one again has
to distinguish two cases: either
$(q_{11}q_{21}^{m'_0-1})^{m'_1}q_{11}^2q_{21}^{-1}=1$ (see Step
4a1) or $q_{11}q_{21}^{m'_0-1}\in R_{m'_1+1}$,
$(q_{11}^2q_{21}^{-1})^{m'_1+1}\not=1$ (see Step 4a2).

Step 4a1. \textit{Assume that $q_{22}=q_{21}^{-1}$,
$q_{11}^6q_{21}^{m'_0-3}=-1$,
$q_{11}^{m'_1+2}q_{21}^{(m'_0-1)m'_1-1}=1$, $q_{11}\in
R_{m'_0+1}$, $q_{21}^{m'_0+1}\not=1$, $m'_0>1$, and $m'_1>1$.}
Then from (\ref{eq-4a-1}) one concludes
\begin{align*}
& & \scp{'}{1}=&1,& \scq{'}{11}{2}=&q_{11},&
\scq{'}{12}{2}\scq{'}{21}{2}=&q_{11}^2q_{21}^{-1},&
\scq{'}{22}{2}=&q_{11}q_{21}^{m'_0-1},\\
m'_2=&m'_0, & \scp{'}{2}=&q_{11}q_{21}^{-1},&
\scq{'}{11}{3}=&q_{21}^{-1},&
\scq{'}{12}{3}\scq{'}{21}{3}=&q_{21},&
\scq{'}{22}{3}=&q_{11},\\
m'_3=&1, & \scp{'}{3}=&1,& \scq{'}{11}{4}=&q_{11},&
\scq{'}{12}{4}\scq{'}{21}{4}=&q_{21},&
\scq{'}{22}{4}=&q_{21}^{-1}.
\end{align*}
Hence the sequence $(m'_i)_{i\in \ndN _0}$ is periodic with period
4. The matrix of the linear map $T'_3T'_2T'_1T'_0$ in Lemma
\ref{l-seq} takes the form
\begin{align}\label{eq-4a-2}
\prod _{i=0}^3\tilde{T}_i\tilde{\tau }=\prod _{i=0}^3
\begin{pmatrix} m'_i & -1 \\ 1 & 0 \end{pmatrix}
=\begin{pmatrix} m'_0{}^2m'_1-m'_0m'_1-2m'_0+1 & m'_0(2-m'_0m'_1)
\\ m'_0m'_1-m'_1-1 & 1-m'_0m'_1 \end{pmatrix}.
\end{align}
By Lemma \ref{l-seq} this matrix has to have finite order. By
Lemma \ref{l-finord} this gives that $0\le (m'_0-2)(m'_0m'_1-2)\le
4$. Since $m'_0\ge 2$ and $m'_1\ge 2$, the latter means that
$m'_0=2$ or $m'_0=3$, $m'_1=2$. In the first case the matrix in
(\ref{eq-4a-2}) has trace $-2$, but it is not equal to $-\id $,
and in the second its trace is $2$, but the matrix is not the
identity. Thus in both cases one obtains a contradiction to Lemma
\ref{l-finord}.

Step 4a2. \textit{Consider the setting $q_{22}=q_{21}^{-1}$,
$q_{11}^6q_{21}^{m'_0-3}=-1$, $q_{11}\in R_{m'_0+1}$,
$q_{11}q_{21}^{m'_0-1}\in R_{m'_1+1}$, $q_{21}^{m'_0+1}\not=1$,
$(q_{11}^2q_{21}^{-1})^{m'_1+1}\not=1$, $m'_0>1$, and $m'_1>1$.}
Since $m'_0>1$ and $m'_1>1$, Lemma \ref{l-no1} implies that $m'_0,
m'_1\in \{2,3\}$.

If $m'_0=2$ then $q_{11}\in R_3$ and hence Equation
$q_{11}^6q_{21}^{m'_0-3}=-1$ implies that $q_{21}=-1$. In this
case the relation $q_{11}q_{21}^{m'_0-1}\in R_{m'_1+1}$ is a
contradiction to $m'_1\in \{2,3\}$.

If $m'_0=3$ and $m'_1=2$ then $q_{11}^2=-1$ and $q_{11}q_{21}^2\in
R_3$. The $6^\mathrm{th}$ power of the latter relation gives that
$q_{21}^{12}=-1$. Thus $q_{21}\in R_8\cup R_{24}$. If $q_{21}\in
R_8$ then $q_{21}^2\in R_4=\{q_{11},-q_{11}\}$ which is a
contradiction to $q_{11}q_{21}^2\in R_3$. On the other hand, if
$q_{21}\in R_{24}$ then $q_{11}\in \{q_{21}^6,-q_{21}^6\}$. Again
$q_{11}=-q_{21}^6$ is a contradiction to $q_{11}q_{21}^2\in R_3$.
The remaining case, where $q_{21}\in R_{24}$ and
$q_{11}=q_{21}^6$, appears in line 13 of Figure 1.

If $m'_0=m'_1=3$ then again $q_{11}^2=-1$. Equation
$q_{11}q_{21}^{m'_0-1}\in R_{m'_1+1}$ gives that
$(q_{11}q_{21}^2)^2=-1$, that is $q_{21}^4=1$. This is a
contradiction to $q_{21}^{m'_0+1}\not=1$.

Step 4b. \textit{Assume now that $q_{22}=-1$, $q_{21}^2\not=1$,
$m'_0>1$, and $m'_1>1$.} There are again two cases:
$q_{11}^{m'_0}q_{21}=1$ or $q_{11}\in R_{m'_0+1}$. In the first
one gets $\scp{'}{0}=1$ and $\scq{'}{11}{1}=-1$ which is a
contradiction $m'_1>1$. Therefore one has $q_{11}\in R_{m'_0+1}$,
and hence
\begin{align}\label{eq-4b-1}
\scp{'}{0}=&q_{11}^{-1}q_{21},&
\scq{'}{11}{1}=&-q_{11}q_{21}^{m'_0},&
\scq{'}{12}{1}\scq{'}{21}{1}=&q_{11}^2q_{21}^{-1},&
\scq{'}{22}{1}=&q_{11}.
\end{align}
{}From Lemma \ref{l-no1} and the above equations one obtains that
$m'_0,m'_1\in \{2,3\}$ and that
$-q_{11}q_{21}^{m'_0}(q_{11}^2q_{21}^{-1})^2q_{11}=-1$. Further,
by Equation (\ref{eq-mij}) one has to have
$(-q_{11}q_{21}^{m'_0})^{m'_1}q_{11}^2q_{21}^{-1}=1$ (see Step
4b1) or $-q_{11}q_{21}^{m'_0}\in R_{m'_1+1}$,
$(q_{11}^2q_{21}^{-1})^{m'_1+1}\not=1$ (see Step 4b2).

Step 4b1. \textit{The setting $q_{11}\in R_{m'_0+1}$, $q_{22}=-1$,
$(-1)^{m'_1}q_{11}^{m'_1+2}q_{21}^{m'_0m'_1-1}=1$,
$q_{11}^6q_{21}^{m'_0-2}=1$, $q_{21}^2\not=1$, $m'_0,m'_1\in
\{2,3\}$.} If $m'_0=3$ then $q_{11}^2=-1$ and $q_{11}^6q_{21}=1$,
which gives a contradiction to $q_{21}^2\not=1$.

If $m'_0=2$ and $m'_1=2$ then $q_{11}\in R_3$ and
$q_{11}q_{21}^3=1$, and hence $q_{21}\in R_9$,
$q_{11}=q_{21}^{-3}$. This example appears in line 10 of Figure 1.

If $m'_0=2$ and $m'_1=3$ then $q_{11}\in R_3$ and
$q_{11}^{-1}q_{21}^5=-1$. Thus $q_{21}^{15}=-1$ and
$q_{21}^5\not=-1$, and hence $q_{21}\in R_{30}\cup R_6$. If
$q_{21}\in R_{30}$ then we have got an example from line 15 in
Figure 1. Otherwise the equations $m''_0=1$, $\scq{''}{11}{0}=-1$,
$\scq{''}{12}{0}\scq{''}{21}{0}=q_{21}$, $\scq{''}{22}{0}=q_{11}$,
and (\ref{eq-ppprime}), (\ref{eq-qpprime}) imply that
$\scp{''}{0}=-q_{21}$ and $\scq{''}{11}{1}=-q_{21}q_{11}=1$,
$\scq{''}{12}{1}\scq{''}{21}{1}=q_{21}^{-1}$. In this case
(\ref{eq-mij}) has no solution for $m''_1$, which is a
contradiction by Lemma \ref{l-seq}.

Step 4b2. \textit{The setting $q_{11}\in R_{m'_0+1}$, $q_{22}=-1$,
$q_{11}^6q_{21}^{m'_0-2}=1$, $-q_{11}q_{21}^{m'_0}\in R_{m'_1+1}$,
$(q_{11}^2q_{21}^{-1})^{m'_1+1}\not=1$, $q_{21}^2\not=1$,
$m'_0,m'_1\in \{2,3\}$.} If $m'_0=3$ then $q_{11}^2=-1$, and
Equation $q_{11}^6q_{21}^{m'_0-2}=1$ is a contradiction to
$q_{21}^2\not=1$.

If $m'_0=2$ and $m'_1=3$ then $q_{11}\in R_3$ and
$-q_{11}q_{21}^2\in R_4$. Therefore $q_{11}^2q_{21}^4=-1$, and
hence $q_{11}=-q_{21}^4$. Using Relation $q_{11}\in R_3$ one gets
$q_{21}^{12}=-1$ and $q_{21}^4\not=-1$, and hence $q_{21}\in
R_{24}$. This example appears in line 13 of Figure 1.

Finally, if $m'_0=m'_1=2$ then $q_{11}\in R_3$ and
$-q_{11}q_{21}^2\in R_3$. Since there are two elements in $R_3$,
there are two possibilities. The first case is when
$q_{11}=-q_{11}q_{21}^2$, that is $q_{21}\in R_4$. Set $\zeta
:=q_{11}^{-1}q_{21}^{-1}\in R_{12}$. Then $q_{11}=-\zeta ^2$ and
$q_{21}=\zeta ^3$, and we obtained an example from line 9 of
Figure 1. On the other hand, if $q_{11}^2=-q_{11}q_{21}^2$ then
$q_{11}=-q_{21}^2$. Since $q_{11}\in R_3$, this gives that
$q_{21}^6=-1$ and $q_{21}^2\not=-1$, that is $q_{21}\in R_{12}$.
This is an example from line 8 of Figure 1.
\end{bew}

\begin{bems}
1. For the if part of the assertion in Theorem \ref{t-class2} one
can avoid to use the complicated proof of the finiteness of $\cB
(V)$ in \cite{a-Heck04a}. Indeed, by Lemma \ref{l-finord} it is
sufficient to compute the periods of the sequences $(T'_i\cdots
T'_2T'_1)_{i\in \ndN _0}$ and $(T''_i\cdots T''_2T''_1)_{i\in \ndN
_0}$ for all $V$ appearing in Figure 1, which is much less
technical.

2. Theorem \ref{t-class2} gives also a new proof of the
classification result presented in \cite{a-Heck04b}. The Nichols
algebra $\cB (V)$ is finite dimensional if and only if $\proots
(\cB (V))$ is finite and all \PBW generators have finite height.
By property (P) the latter happens if and only if
\begin{align*}\tag{$*$}
\mathrm{ord}\,\chi (d,d)<\infty \quad \text{for all }d\in \proots
(\cB (V)).
\end{align*}
Having once the list of Nichols algebras in Figure 1, Theorem
\ref{t-corr} tells how to determine $\proots (\cB (V))$. The
finiteness condition ($*$) is then easily checked.
\end{bems}

\bibliographystyle{mybib}
\bibliography{quantum}

\end{document}